\documentclass{amsart}

\usepackage{amssymb,latexsym} \usepackage{palatino}
\usepackage{mathrsfs} \input{sommers.tex}

\usepackage{fullpage}

\newcommand{\nilrad}{\mathfrak{n}} 
 \newcommand{\g}{\mathfrak{g}}
\newcommand{\bo}{\mathfrak{b}}

\newcommand{\ph}{\Tilde{\Pi}}
\newcommand{\aff}{\Tilde}
\newcommand{\posroots}{\ro^+}
\newcommand{\ideal}{\mathcal I}
\newcommand{\al}{{\alpha}}

\newcommand{\coroot}{{\scriptscriptstyle\vee}}
\newcommand{\ro}{{\Phi}}
\newcommand{\co}{{\Phi^{\scriptscriptstyle\vee}}}

\newcommand{\colat}{{Q^{\scriptscriptstyle\vee}}}
\newcommand{\supp}{\text{ supp}}

\title{$B$-stable ideals in the nilradical of a Borel subalgebra}

\author{Eric N. Sommers} 
\address{University of
Massachusetts---Amherst\\ Amherst, MA 01003}
\address{
Institute for Advanced Study \\
Princeton, NJ 08540}
\email{esommers@math.umass.edu}

\thanks{The author was supported in part by NSF grants DMS-0201826 and DMS-9729992.}
\thanks{The author thanks C. Athanasiadis for pointing out an error
in an earlier version of the paper}

\date{March 14, 2003; September 17,2003}

\begin{document}

\begin{abstract}
We count the number of strictly positive $B$-stable ideals in the
nilradical of a Borel subalgebra and prove that
the minimal roots of any $B$-stable ideal are conjugate
by an element of the Weyl group to a subset of the simple roots.
We also count the number of ideals whose minimal roots are conjugate
to a fixed subset of simple roots.
\end{abstract}

\maketitle

\section{Introduction}

Let $G$ be a connected 
simple algebraic group over the complex numbers and $B$ a Borel
subgroup of $G$.  
Let $\g$ be the Lie algebra of $G$ and $\bo$ the Lie
algebra of $B$.  The nilradical of $\bo$ is denoted $\nilrad$.

The subspaces of $\nilrad$ which are stable under the adjoint action of $B$ are
called $B$-stable ideals of $\nilrad$ (henceforth,
often called ideals).  
The study of these ideals has attracted much recent
attention in the work of Kostant, Cellini-Papi, Panyushev, and others.

The purpose of this note is to extend the recent uniform proof of Cellini-Papi
on the number of ideals to the number of strictly positive ideals (these
are the ones which intersect the simple root spaces trivially).  At the same
time, we obtain a result on the minimal roots in an ideal:  namely,
a set of mutually inequivalent positive roots is conjugate
by an element of the Weyl group to a subset of the simple roots.
We plan to use this result in a later paper to study Kazhdan-Lusztig cells. 
Finally, we count the number of ideals whose minimal roots are conjugate to a
fixed set of simple roots.

\section{Statement of Results}

Fix a  maximal torus $T$ in $B$ and let
$(X, \ro, Y, \co)$ be the root datum determined
by $G$ and $T$, and let $W$ be the Weyl group.
Let $\Pi \subset \posroots$ denote the simple roots and positive roots 
determined by $B$.  As usual, $\langle \ , \ \rangle$ denotes the pairing
of $X$ and $Y$.  Let $\colat$ denote the lattice
in $Y$ generated by $\co$ (the coroot lattice). 
We denote the
standard partial order on $\ro$ by $\prec$; so $\al \prec \beta$ for 
$\al, \beta \in \ro$ if and only if $\beta - \al$ is a sum of positive roots.

We define an ideal (also called an upper order ideal) 
$\ideal$ of $\posroots$ to be a collection
of roots such that if 
$\alpha \in \ideal, \beta \in \posroots$, and $\alpha+\beta \in \posroots$,
then $\alpha + \beta \in \ideal$.
In other words, if $\al \in \ideal$ and $\gamma \in \posroots$ with $\al \prec \gamma$,
then $\gamma \in \ideal$. 

It is easy to see that $B$-stable ideals in the nilradical $\nilrad$ of $\bo$
are naturally in bijection with the ideals of $\posroots$.
Namely, if $I$ is a $B$-stable ideal of $\nilrad$, it is stable under the action of $T$,
hence $I$ is a sum of roots spaces.  Denote by $\ideal$ the set of roots 
whose root space is contained in $I$.  Then $\ideal$ is an ideal of $\posroots$ 
and this map is a bijection.      

Given an ideal $\ideal$ in $\posroots$, we define the minimal roots $\ideal_{min}$ 
of $\ideal$ as follows:  $\al \in \ideal$ belongs to $\ideal_{min}$ 
if and only if $\beta \in \posroots, \beta \prec \al$ implies $\beta \notin \ideal$.
Clearly $\ideal$ determines and is determined by its set of minimal roots $\ideal_{min}$. 
Note that the elements of $\ideal_{min}$ are mutually inequivalent
elements of $\posroots$ and that every set of mutually inequivalent elements
is an $\ideal_{min}$ for a unique $\ideal$ (namely, $\ideal$
is the set of all elements bigger or equal to the elements of $\ideal_{min}$).

Our two main results are the following:

\begin{thm2} 
Let $\ideal$ be an ideal of $\posroots$.
Then there exists $w \in W$ such that $w(\ideal_{min}) \subset \Pi$.
\end{thm2}

In other words, any set of mutually inequivalent elements of $\ro^+$ 
is conjugate by an element of $W$ to a subset of $\Pi$.  

We say an ideal $\ideal$ of $\posroots$ is {\bf strictly positive} 
if $\ideal \cap \Pi$ is empty.

\begin{thm2}
The number of strictly positive ideals is given by
$$\frac{1}{|W|} \prod_{i=1}^{n} (h - 1 +m_i)$$
where $m_1, \dots, m_n$ are the exponents of $W$
and $h$ is the Coxeter number.
\end{thm2}

Cellini-Papi gave a uniform proof of the corresponding result for all ideals (where $h-1$
in the formula gets replaced by $h+1$)  \cite{cellini-papi:1}, \cite{cellini-papi:2}.
Earlier authors had counted the ideals, but had not produced this closed formula
(see \cite{shi:2}). 
We will review the work of \cite{cellini-papi:1}, \cite{cellini-papi:2}, and also
that of \cite{shi:1},
since it will be needed for the theorem on minimal roots.  

The formula for the number of strictly positive ideals 
shows up in the work of Fomin and Zelevinsky \cite{fomin-zel:1},
where it counts the number of postive clusters.  That work was an inspiration for the present one.

Finally, 
we note that Athanasiadis has obtained our second theorem (by a different method)
\cite{christos}, 
as has Panyushev, who has also obtained some of our results on minimal roots
\cite{panyushev:2}.

\section{Roots lemmas}

We first prove a few lemmas about roots.

\begin{lem} \label{first_lemma}
Let $\gamma \in \ro$ and suppose 
$\gamma = \sum_{i=1}^k \al_i$ for $\al_i \in \ro$.  For $1 \leq j \leq k$,
either $\gamma - \al_j \in \ro \cup \{ 0 \}$ or $\al_j+\al_l \in \ro \cup \{ 0 \}$ for some $l$ with 
$1 \leq l \leq k$, $l \neq j$.
If in the latter case $l$ is unique and $\al_j \neq -\al_l$, then 
$\al_l$ is a long root and $\al_j$ is short (and the root system is not simply-laced). 
\end{lem}

\begin{proof}
If $\langle \gamma, \al_j^{\scriptscriptstyle\vee} \rangle > 0$, then applying the reflection
$s_{\al_j}$ to $\gamma$
ensures that $\gamma - \al_j \in \ro $ (since strings of roots 
are unbroken), unless $\gamma =\al_j$ in which case $\gamma - \al_j = 0$.
Therefore 
$\gamma - \al_j \notin \ro \cup \{ 0 \}$ implies 
$\langle \gamma, \al_j^{\scriptscriptstyle\vee} \rangle \leq 0$.
This implies that $\langle \al_l , \al_j^{\scriptscriptstyle\vee} \rangle < 0$ for some $l \neq j$
since $\langle \al_j , \al_j^{\scriptscriptstyle\vee} \rangle = 2$.
And this in turn forces $\al_l + \al_j \in \ro \cup \{ 0 \}$ as in the first sentence of the proof
with $\al_j$ replaced by $-\al_j$.
Finally if $\al_l$ is short or $\al_j$ is long and $\al_j \neq -\al_l$, then 
$\langle \al_l , \al_j^{\scriptscriptstyle\vee} \rangle = -1$; hence there exist $m \neq l$
with $\langle \al_m , \al_j^{\scriptscriptstyle\vee}  \rangle < 0$, implying
that also $\al_m + \al_j \in \ro \cup \{ 0 \}$.
\end{proof}



\begin{lem} \label{reorder}
Let $\al_1 \in \ro \cup \{ 0 \}$.  Suppose that 
$\al_i \in \posroots$ for $2 \leq i \leq k$ and $\gamma = \sum_{i=1}^k \al_i \in \ro \cup \{ 0 \}$.
Then there exists a re-ordering of the $\al_i$'s with $i \geq 2$ so that 
$\sum_{i=1}^j \al_i \in \ro \cup \{ 0 \}$ for all $j \in \{ 1, \dots, k\}$.
\end{lem}

\begin{proof}
We proceed by induction on $k$, the case $k=2$ being trivial.
Assume $k \geq 3$.  If $\gamma - \al_k \in \ro \cup \{ 0 \}$, we are done by induction
applied to $\gamma - \al_k$, which is a sum of $k-1$ roots.
If not, then $\al_k + \al_l \in \ro \cup \{ 0 \}$ 
for some $l$ by Lemma \ref{first_lemma}.
If $l=1$, we finish by applying induction to the sum of the $k-1$ roots
$(\al_1 + \al_k) + \sum_{i \neq 1, k} \al_i = \gamma$
and then breaking apart the first two roots.
If $l>1$, then $\al_k + \al_l \in \posroots$ and 
we apply induction to the sum of the $k-1$ roots
$\al_1 + (\al_k+\al_l) + \sum_{i \neq 1, k, l } \al_i = \gamma$.
At some point in the re-ordering of the latter $k-2$ roots, we find that
$\beta$, $\beta + (\al_k+\al_l) \in \ro \cup \{ 0 \}$ 
where $\beta$ is a sum of $\al_i$'s.  
We apply the previous lemma with $\gamma = \beta + \al_k+ \al_l$.  
Now $\al_k \neq -\al_l$ since $\al_k$ and $\al_l$ are positive roots.
Thus if neither $\beta + \al_k$ nor $\beta + \al_l$ is in 
$\ro \cup \{ 0 \}$, the previous lemma would imply that $\al_k$ is both long and short
(in a root system with distinct lengths), a contradiction.
Therefore either $\beta + \al_k$ or $\beta + \al_l$ is in 
$\ro \cup \{ 0 \}$, completing the proof.
\end{proof}

To each $\ideal$ and each  $\al \in \posroots$, we now attach a nonnegative integer
(two integers if $\ideal$ is strictly positive).
One can use the previous lemma to show that the numbers
in the first part of the definition below coincide with the numbers attached to
each root and each ideal in \cite{cellini-papi:1}. 

\begin{defn} \label{define_min}
Let $\ideal$ be an ideal and $\al \in \posroots$.
\begin{enumerate}
\item Let $\al_{\ideal, +}$ be defined as follows
$$\al_{\ideal, +}:= \text{max} \ \{ k \ | \ \al = \sum_{i=1}^k \gamma_i \text{ with }
\gamma_i \in \ideal \}.$$
\item If $\ideal$ is strictly positive, let $\al_{\ideal, -}$ be defined as follows
$$\al_{\ideal, -}:= \text{min} \ \{ k \ | \ \al = \sum_{i=1}^{k+1} \gamma_i 
 \text{ with } \gamma_i \in \posroots \! - \! \ideal \}.$$
\end{enumerate}
\end{defn}

The fact that $\ideal$ is strictly positive ensures that
$\al_{\ideal, -}$ is defined since $\Pi \subset \posroots \! - \! \ideal$.    

\section{Affine Weyl group}

Let $W_a$ be the affine Weyl group of $G$.  
Then $W_a$ has two descriptions.  
First, it is isomorphic to $W \ltimes \colat$
and it acts naturally on 
the vector space $V = \colat \otimes \mathbf{R}$, with 
$W$ acting in the usual way and $\colat$ acting by translations.
If $\lambda \in \colat$, we write $\tau_{\lambda}$ for 
the corresponding element of $W_a$.
Second, $W_a$ is a Coxeter group.  Let $S = \{ s_0, s_1, \dots, s_n \}$
be the simple reflections of $W_a$, where $s_0$ is the
affine reflection $s_{\theta} \tau_{-\theta^\coroot}$
and $s_i$ for $i >0$ is the reflection via the simple root $\al_i \in \Pi$.
Here, $\theta$ is the highest root of $\posroots$.
$W_a$ comes equipped with a length function $l( -)$ and the Bruhat order.

The affine roots of $W_a$ are denoted $\al + m \delta$
where $\al \in \ro$ and $m$ is an integer and $\delta$ is, for our purposes,
just a place-keeper for $m$.
The positive affine roots are 
those $\al + m \delta$ with either $\al \in \posroots$ and $m>0$ 
or $-\al \in \posroots$ and $m \geq 0$.
If $\aff{\al}$ is an affine root, we write $\aff{\al} \succ 0$
if $\aff{\al}$ is positive and $\aff{\al} \prec 0$ otherwise.
The simple affine roots are 
those $-\al$ with $\al \in \Pi$ together with $\theta + \delta$.
If $w \in W_a$ is written $w_f \tau_{\lambda}$,
then the action of $w$ on the affine roots is given by
$$w (\al + m \delta) = w_f(\al) + (m + \langle \al, \lambda \rangle) \delta.$$
This action is consistent with the action of $W_a$ on $V$ as follows.
Consider the hyperplanes 
$$H_{\al, m} = \{ v \in V \ | \ \langle \al, v \rangle = m \}.$$
Then $w(H_{\al, m}) = H_{\beta, k}$ whenever 
$w(\al + m \delta) = \beta + k \delta$.  

Let 
$$A = \{ v \in V \ | \ \al(v) > 0 \text{ for } \al \in \Pi \text{ and } \theta(v) < 1 \}.$$
Given $w \in W_a$, let $N(w)$ denote the positive
affine roots $\aff{\al} = \alpha + m \delta$ such that $w^{-1}(\aff{\al}) \prec 0$.
This is equivalent to the following: 
$\al + m \delta \in N(w)$ if and only if $H_{\al, m}$ 
separates $A$ from $w(A)$. 
We define the support of $w$, denoted $\supp (w)$, 
to be those roots $\al \in \ro$ 
such that $\al + m \delta \in N(w)$ for some $m$.
We say $w$ is {\bf dominant} if $\supp(w) \subset \posroots$.  This
is equivalent to $w(A)$ belonging to the dominant chamber 
$C := \{ v \in V \ | \ \al(v) > 0 \text{ for } \al \in \Pi \}$.

For a point $v \in V$ which is not on any affine hyperplane $H_{\beta, m}$, define
$k(\al, v) \in \mathbf{Z}$ for $\al \in \posroots$ so that 
$k(\al, v) < \langle \al, v \rangle < k(\al, v) + 1$.
For $w \in W_a$ define
$k(\al, w) =  k(\al, v)$ for any point $v$ in $w(A)$. 
The following result is important and goes back to Shi \cite{shi:1} (see also the
references in \cite{cellini-papi:1} to their earlier work).  Shi's formula looks
a bit more complicated since he uses coroots instead of roots.  

\begin{prop} \label{keyaffine} \cite{shi:1}
Let $k_{\al}$ be a collection of integers for each $\al \in \posroots$.
Then there exists $w \in W_a$ with $k(\al, w) = k_{\al}$ for all $\al \in \posroots$
if and only if whenever $\al, \beta, \al +\beta \in \posroots$ 
the inequalities 
$$k_{\al} + k_{\beta} \leq k_{\al + \beta} \leq k_{\al} + k_{\beta} +1$$
hold true.
\end{prop}

We will use a couple of standard facts about $W_a$.
First, the length of $l(w)$ of $w \in W_a$ is the cardinality of $N(w)$.
Second, $l(w) = \sum |k(\al, w)|$. 
Third, $N(x) \subset N(y)$ if and only if $y = x u$ 
for some $u \in W_a$ with $l(y) = l(x) + l(u)$. 
Fourth, $ws_i < w$ in the Bruhat order 
if and only if $w(\aff{\al}_i)$ is negative where $s_i$ is the reflection
corresponding to the affine simple root $\aff{\al_i}$.
We also note that $W_a$ acts simply-transitively on the set of all alcoves
(that is, the regions in $V$ of the form $w(A)$ for $w \in W_a$).


\medskip

Given an ideal $\ideal$ in $\posroots$, let $ST(\ideal) \subset W_a$
denote all $w$ such that $\supp (w)=\ideal$. 
The $ST$ refers to sign type (after Shi). 
It follows from Lemma \ref{what_is_min} 
that $ST(\ideal)$ is always non-empty.  

Recall that a strictly positive ideal is one for which $\ideal \cap \Pi$ is empty.
The strictly positive ideals are exactly the ideals for 
which $ST(\ideal)$ is a finite set.  
Namely, if $\al \in \ideal \cap \Pi$,
then $ST(\ideal)$ contains the elements $\tau_{m \omega^{\vee}} w$
for $m \geq 0$ where $\omega^{\vee}$ is the fundamental coweight corresponding to 
$\al$ and $w$ is any element of  $ST(\ideal)$ (we require that $m \omega^{\vee} \in
\colat$). Hence  $ST(\ideal)$ is not finite.  
Conversely, if $\ideal \cap \Pi$ is empty and 
$w \in ST(\ideal)$, then $w(A)$ lies in the region
bounded by the hyperplanes $\{ H_{\al,0} \cup  H_{\al,1} \ | \ \al \in 
\Pi \}$, a region
of finite volume. 

\section{Enumerating and counting ideals}

\begin{lem} \label{max/min}
Let $\ideal$ be an ideal in $\posroots$.
\begin{enumerate}
\item \cite{cellini-papi:1} $w \in ST(\ideal)$ implies $k(\al, w) \geq 
\al_{\ideal, +}$ for all $\al \in \posroots$.
\item If $\ideal$ is strictly positive then 
$w \in ST(\ideal)$ implies $k(\al, w) \leq 
\al_{\ideal, -}$ for all $\al \in \posroots$.
\end{enumerate}
\end{lem}

\begin{proof}
If $\gamma \in \ideal$, then certainly $k(\gamma,w) \geq 1$.
So if $\al =  \sum_{i=1}^k \gamma_i$ with $\gamma_i \in \ideal$,
then by Proposition \ref{keyaffine}, $k(\al,w) \geq k$.  Hence $k(\al, w) \geq 
\al_{\ideal, +}$.  

Similarly, if $\ideal$ is strictly positive
and if $\gamma \in \posroots - \ideal$, then 
$k(\gamma,w)=0$.  So if $\al =  \sum_{i=1}^k \gamma_i$ with $\gamma_i \in 
\posroots - \ideal$, then 
$k(\al,w) \leq k-1$ by Proposition \ref{keyaffine}.  Therefore, $k(\al, w) \leq 
\al_{\ideal, -}$.
\end{proof}

\begin{lem} \label{what_is_min}
Let $\ideal$ be an ideal.
\begin{enumerate}
\item \cite{shi:1}, \cite{cellini-papi:1} 
There exists $w \in W_a$ such that $k(\al, w)=
\al_{\ideal, +}$ for all $\al \in \posroots$.
\item If $\ideal$ is strictly positive then 
there exists $w \in W_a$ such that $k(\al, w)=
\al_{\ideal, -}$ for all $\al \in \posroots$.
\end{enumerate}
\end{lem}

\begin{proof}
Suppose $\al, \beta, \al+\beta \in \posroots$.

For the first statement, let $a=\al_{\ideal, +}, b= \beta_{\ideal, +},
c= (\al+\beta)_{\ideal, +}$ for simplicity.
By Proposition \ref{keyaffine} we need to show that $a+b \leq c \leq a+b+1$.
The first inequality is obvious from the definitions.  For the second inequality,  
write $\al + \beta =  \sum_{i=1}^k \gamma_i$ where $\gamma_i
\in \ideal$.  Then 
$\al = -\beta + \sum_{i=1}^k \gamma_i$.
By Lemma \ref{reorder} there exists $j$ (after re-ordering the $\gamma_i$
as in the lemma)
for which $-\beta + \sum_{i=1}^{j-1} \gamma_i \in \ro^- \cup \{ 0 \}$ and 
$-\beta + \sum_{i=1}^j \gamma_i \in \ro^+$.
The former inclusion means that $\mu + \sum_{i=1}^{j-1} \gamma_i = \beta$
for some $\mu \in \posroots \cup \{ 0 \}$.  Hence $b \geq j-1$ 
since $\mu + \gamma_i \in \posroots$ for some $i \leq j-1$ by Lemma \ref{reorder}
and then $\mu + \gamma_i \in \ideal$ as $\gamma_i \in \ideal$ and $\mu \in \posroots \cup \{ 0 \}$.
The inclusion $-\beta + \sum_{i=1}^j \gamma_i \in \ro^+$ implies
$-\beta + \sum_{i=1}^{j+1} \gamma_i \in \ideal$ since $\gamma_{j+1} \in \ideal$
(unless of course $j=k$).
Thus $a \geq k-j$ (which is also true if $j=k$). 
Together, $a+b \geq k-1$, so $a+b+1 \geq c$.

For the second statement, let $a=\al_{\ideal, -}, b= \beta_{\ideal, -},
c= (\al+\beta)_{\ideal, -}$.
Clearly from the definitions, $c \leq a+b+1$.
For the other inequality, 
write $\al + \beta =  \sum_{i=1}^k \gamma_i$ for $\gamma_i
\in \posroots - \ideal$,
and let $j$ be as in the first part of the proof.
The definition of $j$ implies that 
$-\beta + \sum_{i=1}^j \gamma_i 
= (-\beta + \sum_{i=1}^{j-1} \gamma_i) + \gamma_j \in \ro^+ - \ideal$
since the term in parentheses is either zero or a negative root and $\gamma_j \in \ro^+ - \ideal$.
Thus $a \leq k-j$.
On the other hand, $-\beta + \sum_{i=1}^j \gamma_i \in \ro^+$ means that
$\beta = -\mu + \sum_{i=1}^j \gamma_i$ for some $\mu \in \posroots$.
By Lemma \ref{reorder} applied to $\beta$,
we can re-arrange the $\gamma_i$'s with $i \leq j$ so that for some $l \leq j$ we have
$-\mu + \sum_{i=1}^{l-1} \gamma_i \in \ro^- \cup \{ 0 \}$ and
$-\mu + \sum_{i=1}^{l} \gamma_i \in \posroots$.
It follows that $-\mu + \sum_{i=1}^l \gamma_i \in \posroots - \ideal$
since $\gamma_l \in \posroots - \ideal$.  We conclude which that $b \leq j-1$.
Together, so $a+b \leq k-1$ and hence
$a+b \leq c$. 
\end{proof}

Let $w_{\ideal, min}$ denote the unique element of minimal length in $ST(\ideal)$
and when $\ideal$ is strictly positive let 
$w_{\ideal, max}$ denote the unique element of maximal length in $ST(\ideal)$.
These exist and are unique by the previous two lemmas (and the simple-transitivity
of the affine Weyl group on alcoves).  We drop the $\ideal$ from the subscript
when there is no confusion.

\begin{rmk}
Let $R$ be a region of the Shi arrangement.  Shi showed more generally 
that there is a unique element $w$ with $w(A) \subset R$ such that, 
for every positive root $\al$, $|k(\al,w)| \leq |k(\al, w')|$ 
for every $w'$ with $w'(A) \subset R$.   
Cellini-Papi characterized these minimal numbers $k(\al, w)$ when
$w(A)$ lies in the dominant chamber, allowing for an alternative (and simpler) proof
in that case.

It is also possible to study (in the same spirit as Shi) the maximal elements of all
bounded regions, which we have only done here (in the spirit of Cellini-Papi) for those
bounded regions in the dominant chamber (that is, those ideals $\ideal$ 
for which $ST(\ideal)$ is finite).
\end{rmk}

\begin{prop} \label{key_shi}
Let $\ideal$ be an ideal and $w \in ST(\ideal)$.
\begin{enumerate}
\item \cite{shi:1} 
$w = w_{min}$ if and only if for all $s \in S$ for which $ws<w$ we have 
$ws \notin ST(\ideal)$
\item 
$\ideal$ is strictly positive and $w \! = \! w_{max}$ 
if and only if for all $s \in S$ for which $ws>w$ we have $ws \notin ST(\ideal)$
\end{enumerate}
\end{prop}

\begin{proof}

For the forward direction of the first statement, assume $w = w_{min}$.
Then $w s < w$ implies 
$k(w s, \beta) = k(w, \beta) - 1$ for a unique $\beta \in \posroots$
by the basic properties of $W_a$ listed above.  This implies that
$w s  \notin ST(\ideal)$ by the previous two lemmas
(and so, in fact, $k(w s, \beta) = 0$ and $\beta \in \ideal_{min}$).

For the reverse direction of the first statement, assume $w \neq w_{min}$.
Then since $w \in ST(\ideal)$, by Lemma \ref{max/min} we have $N(w_{min}) \subset N(w)$.
Thus by the third property of the affine Weyl group listed above, 
there exists $u \in W_a$, $u \neq 1$ so that $w = w_{min} u$ and 
$l(w) = l(w_{min}) + l(u)$.  It follows that 
there exists $s \in S$ so that $ws < w$ and $ws$ has a reduced expression
beginning with $w_{min}$.  Hence $N(w_{min}) \subset N(ws) \subset N(w)$ 
which implies that $ws \in ST(\ideal)$ and the reverse direction is proved.

For the forward direction of the second statement, set $w = w_{max}$ for $\ideal$.
Then $w s > w$ implies 
$k(w s, \al) = k(w, \al) + 1$ for a unique $\al \in \posroots$,
which implies that
$w s  \notin ST(\ideal)$ by the previous two lemmas.

For the reverse direction of the second statement,
we first need to show that $\ideal$ is strictly positive.
Let $w$ satisfy the hypothesis that for all 
$s \in S$ for which $ws>w$ we have $ws \notin ST(\ideal)$.
If $\ideal$ were not strictly positive, then the proof that $ST(\ideal)$
is not finite shows that there exists $x \in ST(\ideal)$ (in fact, infinitely many) 
with $N(w) \subset N(x)$.
Hence there exists $u \in W_a$, $u \neq 1$ so that $x = w u$ and 
$l(x) = l(w) + l(u)$.  Thus there exists $s \in S$ with $w < ws < x$. 
In particular
$N(w) \subset N(ws) \subset N(x)$.  
Therefore $ws \in ST(\ideal)$, contradicting our hypothesis on $w$.
We conclude that $\ideal$ is strictly positive and possesses a maximal element.
 
Now assume that $w \in ST(\ideal)$ but $w \neq w_{max}$.  
By Lemma \ref{max/min}, we have $N(w) \subset N(w_{max})$.
Then there exists $u \in W_a$, $u \neq 1$ so that $w_{max} = w u$ and 
$l(w_{max}) = l(w) + l(u)$.  Hence there exists $s \in S$ with $w < ws < w_{max}$.  
In particular
$N(w) \subset N(ws) \subset N(w_{max})$.  
Therefore $ws \in ST(\ideal)$ finishing the proof of the reverse direction.
\end{proof}

Now let $t$ be a natural number which is relatively 
prime to the Coxeter number $h$ of $G$ and write $t = ah +b$
where $1 \leq b < h$.
Let $\ro_{k}$ be the roots of $\ro$ of height $k$ (the sum of 
the coefficients when expressing a root in the simple root basis).
 
Define 
$$D^t = \{ \lambda \in \colat \ | \ \langle \al, \lambda \rangle \leq a
\text{ for } \al \in \ro_{b} \text{ and }
\langle \al, \lambda \rangle \leq a+1 \text{ for } \al \in \ro_{b-h}\}.$$

We are interested in the two cases where $t= h+1$ and $t=h-1$.
When $t=h+1$, we have $\ro_b = \Pi$ and $\ro_{b-h}= \{ -\theta \}$.
When $t=h-1$, we have $\ro_b = \{ \theta \}$ and $\ro_{b-h}= -\Pi$.

The following is the main result of \cite{cellini-papi:2}.  We give a proof
here in order to extract information on the minimal roots of an ideal.

\begin{prop} \label{biject_1} \cite{cellini-papi:2}
The set of ideals is in bijection with the elements of $D^{h+1}$.  
\end{prop}

\begin{proof}

For $w \in W_a$, we write $w = x \tau_{\lambda}$ with $x \in W$.
Our aim is to show that if $w=w_{min}$ for some ideal $\ideal$, then 
$\lambda \in D^{h+1}$ and conversely, if $\lambda \in D^{h+1}$
then $w= x \tau_{\lambda}$ is a minimal element of an ideal
for the unique $x \in W$ which makes $x \tau_{\lambda}$ dominant.  
By the uniqueness of $x$ and the uniqueness of the minimal element of an
ideal, this will establish the bijection.

Let $\aff{\al}_i$ be an affine simple root and write
$\aff{\al}_i = \al_i + m_i \delta$ where
$\al_i$ is a negative finite simple root 
(respectively, $\theta$) and $m_i=0$ (respectively, $m_i=1$).
Then 
$$w(\aff{\al}_i) = x(\al_i) + (m_i + \langle \al_i, \lambda \rangle) \delta.$$

First assume that $w=w_{min}$ for some ideal $\ideal$.
If $ws_i > w$, then $w (\aff{\al}_i)$ is positive and so certainly 
$\langle \al_i, \lambda \rangle \geq -m_i$.
If $ws_i < w$, then $w(\aff{\al}_i)$ is negative.  
Since $w$ is dominant and $\supp(ws_i) \neq
\supp(w)$ by the previous proposition, 
we must have $w(\aff{\al}_i) = -\beta - \delta$
where $\beta \in \posroots$ and also  $\beta \in \ideal_{min}$.
Consequently, $x(\al_i) = -\beta$ and 
$\langle \al_i, \lambda \rangle = -1 -m_i$.
We conclude in both cases that if $w= w_{min}$ then 
$\langle \al_i, \lambda \rangle \geq -1-m_i$, that is, $\lambda \in D^{h+1}$.
Moreover, when equality holds, $-x(\al_i) \in \ideal_{min}$.

Conversely, assume that $\lambda \in D^{h+1}$.  Then there is a unique $x \in W$
such that $w = x \tau_{\lambda}$ is dominant.  It is characterized by the fact that
$$\langle \beta, \lambda \rangle \geq 0 \text{ if and only if } x(\beta) \succ 0 \text{ and }$$
$$\langle \beta, \lambda \rangle < 0 \text{ if and only if } x(\beta) \prec 0$$
for all $\beta \in \posroots$. 

Now if $\langle \al_i, \lambda \rangle = -1-m_i$, then 
$x(\al_i) \prec 0$ and
certainly $w(\aff{\al}_i) = x(\al_i) -\delta$ is negative.
Thus $ws_i < w$ and also $\supp(ws_i) \cup \{ -x(\al_i) \} = \supp(w)$,
so $\supp(ws_i) \neq \supp(w)$.
If $\langle \al_i, \lambda \rangle = -m_i$, then $x(\al_i) \prec 0$
and $w(\aff{\al}_i) = x(\al_i)$ is a positive affine root.  Thus $ws_i > w$.
Finally if $\langle \al_i, \lambda \rangle > -m_i$, then certainly
$w(\aff{\al}_i)$ is positive and $ws_i > w$. 
This shows that $w$ satisfies the hypotheses of the previous proposition and so
$w$ equals $w_{\ideal, min}$ for some $\ideal$.

This establishes the Cellini-Papi bijection between ideals (and their minimal elements)
and elements of $D^{h+1}$.
\end{proof}

\begin{prop} \label{biject_2}
The set of strictly positive ideals is in bijection with the elements of $D^{h-1}$.
\end{prop}

\begin{proof}
Suppose $w = w_{max}$ for some strictly positive ideal.  
If $ws_i < w$, then $w(\aff{\al}_i)$ is negative and certainly
$\langle \al_i, \lambda \rangle \leq -m_i$.
If $ws_i > w$, then $w(\aff{\al}_i)$ is positive.  Since $w$ is dominant and 
$\supp(ws_i) \neq \supp(w)$ by Proposition \ref{key_shi}, we have either
$w(\aff{\al}_i) = \beta + \delta$ for $\beta \in \posroots$
or $w(\aff{\al}_i)$ is a negative finite simple root.
In the first case, $\langle \al_i, \lambda \rangle = 1-m_i$  and in the second
case, $\langle \al_i, \lambda \rangle = -m_i$.
In all cases then, 
$\langle \al_i, \lambda \rangle \leq 1-m_i$, that is, $\lambda \in D^{h-1}$.
Moreover, if equality holds, then $x(\al_i)$
is a maximal root in $\posroots - \ideal$.

Conversely, suppose $\lambda \in  D^{h-1}$ and 
let $x \in W$ be the unique element
such that $w = x \tau_{\lambda}$ is dominant. 
If $\langle \al_i, \lambda \rangle = 1-m_i$, then $x(\al_i) \succ 0$
and so $ws_i > w$ and also $\supp(ws_i) \cup \{ x(\al_i) \} = \supp(w)$.
If $\langle \al_i, \lambda \rangle = -m_i$, then $x(\al_i) \prec 0$
and so $w ({\aff{\al}_i}) = x(\al_i)$ is a positive affine root and 
so $ws_i > w$.  Since $w$ is dominant, it is clear that 
$\supp(ws_i) \neq \supp(w)$.
Finally if $\langle \al_i, \lambda \rangle < -m_i$, then 
$w{\aff{\al}_i}$ is clearly negative and $ws_i < w$.
This shows that $w$ satisfies the hypotheses of Proposition \ref{key_shi} and so
$w$ equals $w_{\ideal, max}$ for some strictly positive $\ideal$.

This establishes the bijection between strictly positive ideals 
(and their maximal elements)
and elements of $D^{h-1}$.
\end{proof}

It is now possible to enumerate the number of ideals and strictly positive ideals.
The former was done uniformly in \cite{cellini-papi:2}.

\begin{thm} \label{counting1}
The number of ideals (strictly positive ideals) is given by
$$\frac{1}{|W|} \prod_{i=1}^{n} (t+m_i)$$ where 
$t = h+1$ (respectively, $t=h-1$)
and the $m_i$ are the exponents of $W$.
\end{thm}

\begin{proof}
The cardinality of $D^t$ can be computed in general when $t$ is good for $G$. 
Namely, let $\Delta^t$ denote the simplex in $V$ bounded by the hyperplanes 
$\{ H_{\al, a} \ | \ \al \in \ro_b \} \cup \{ H_{\al, a+1} \ | \ \al \in \ro_{b-h} \}$. 
It was stated (without proof) in \cite{sommers:affineweylgp}
that there exists $\tilde{w} \in W_a$ such that  
$\tilde{w}(\Delta^t) = t \bar{A}$ where $\bar{A}$ is the closure of the fundamental alcove
$A$. 
The existence of this element is proved as follows: 
by \cite{fan} (the end of section 2.3)
there is an element $\tilde{w}'$ 
in the extended affine Weyl group with this property and thus 
(for example by 
\cite{cellini-papi:2}, Lemma 1) there is an element 
$\tilde{w} \in W_a$ with the same property.

It follows that $$\tilde{w}(D^t) =  \colat \cap t \bar{A}.$$  
The latter is known to parametrize the orbits of $W$ on $\colat / t \colat$, and
the number of orbits is given by 
$$\frac{1}{|W|} \prod_{i=1}^{n} (t+m_i),$$ 
(see for example \cite{haiman}).
\end{proof}

\begin{rmk}
It is possible to show a more general result:  namely, let $J \subset \Pi$ and
let $W_J$ be the corresponding parabolic subgroup of $W$.  
Then the number of regions in the Shi arrangement lying in a fundamental domain
for the action of $W_J$ is 
$$\frac{1}{|W_J|} \prod_{i=1}^{j} (h+1+m_i) (h+1)^{n-j},$$ 
where the $m_i$ are the exponents of $W_J$ and $n$ is the rank of $G$.
A similar statement holds for bounded regions by replacing $h+1$ with $h-1$.

This recovers under one rubric Shi's original result on the number of regions
($J = \emptyset$, $t=h+1$),
Cellini-Papi's counting of the dominant regions
($J = \Pi$, $t=h+1$), and Headley's counting
of all bounded regions ($J = \emptyset$, $t=h-1$) \cite{headley}. 
This formula appears in \cite{sommers:affineweylgp} as the Euler characteristic
of a partial affine Springer fiber.  The proofs there, together with Shi's work
in \cite{shi:1}, are sufficient to prove the above formula on 
the total number of regions in a $W_J$-fundamental domain.
\end{rmk}

\section{Combinatorics of the minimal roots of an ideal}

We are now able to obtain some new results on the minimal roots of an ideal
by using the work of the previous section.

\begin{defn}
Let $\lambda \in D^{t}$.  Write $t = ah +b$ as above.  Define 
$$B_{t, \lambda} = \{ \al \in \ro_b \ | \ \langle \al , \lambda \rangle = a \}
\cup \{ \al \in \ro_{b-h} \ | \ \langle \al , \lambda \rangle = a+1 \}.$$
\end{defn} 

In other words, $B_{t, \lambda}$ records the root hyperplanes bounding the simplex
defined by $D^{t+1}$ on which $\lambda$ lies.  

\begin{prop} \label{simples}
Let $\ideal$ be an ideal.  Let $w=w_{min}$ for $\ideal$ 
in the first statement below.  If $\ideal$ is strictly positive,
let $w=w_{max}$ for $\ideal$ in the second statement below.
We write $w = x \tau_{\lambda}$ where $x \in W$.
Let $\al \in \posroots$.
\begin{enumerate}
\item 
Then $\al \in \ideal_{min}$ if and only if 
$x^{-1}(\al) \in B_{h+1, \lambda}$.
\item 
Assume $\ideal$ is strictly positive.
Then $\al$ is a maximal root of $\posroots - \ideal$ 
if and only if 
$x^{-1}(\al) \in B_{h-1, \lambda}$.
\end{enumerate}
\end{prop}

\begin{proof}
We showed in the proofs of Propositions \ref{biject_1}, \ref{biject_2}
that the reverse implication 
in both cases holds.  We now prove the forward implications.

For the first statement,
suppose $\al$ is minimal in $\ideal$.  Then certainly $k(\al, w)=1$.  
We define integers $k_\beta$ for $\beta \in \posroots$ as
follows:  $k_{\beta} = k(\beta, w)$ for $\beta \in \posroots - \al$ and
$k_{\al}=0$.  We wish to show that there exists an element $y \in W_a$ with
$k(\beta,y) = k_{\beta}$ for all $\beta \in \posroots$.
By Proposition \ref{keyaffine} we have to show that if $\gamma \in \posroots$
and $\al + \gamma \in \posroots$, then 
$k_{\al} + k_{\gamma} \leq k_{\al + \gamma} 
\leq  k_{\al} + k_{\gamma} +1$.
In other words, since $k_{\al}=0$, we need 
$k_{\gamma} \leq k_{\al + \gamma} 
 \leq  k_{\gamma} +1$.  
The first inequality is clear (and is clearly a strict inequality).
For the second inequality, 
write $\al + \gamma = 
\sum_{i=1}^k \gamma_i$ with $\gamma_i \in \ideal$.
Then $\gamma = -\al + \sum_{i=1}^k \gamma_i$.  Invoking Lemma \ref{reorder}, 
we have $-\al + \gamma_j \in \ro \cup \{ 0 \}$ for some $j$.  But $\al$ is minimal in $\ideal$,
so $-\al + \gamma_j \in \posroots \cup \{ 0 \}$; hence pulling out one more root
gets us $-\al + \gamma_j + \gamma_{j'} \in \ideal$.  
Thus $k(\gamma, w) \geq  k-1$ and so $k(\gamma,w) \geq k(\al+\gamma, w) -1$.
In other words, $k_{\gamma}+1 \geq k_{\al + \gamma}$, as desired.

It follows that there exists $y \in W_a$ such that $w = y s_i$ 
where $s_i \in S$.  Moreover, $w(\aff{\al}_i) = -\al - \delta$.
Hence $x(-\al_i) = \al$ and $-\al_i \in  B_{h+1, \lambda}$.

For the second statement, suppose $\al$ is maximal in $\posroots - \ideal$.  
Let $k_{\beta} = k(\beta, w)$ for $\beta \in \posroots- \al$ and
$k_{\al}=1$.  We wish to show that there exists an element $y \in W_a$ with
$k(\beta,y) = k_{\beta}$ for all $\beta \in \posroots$.
By Proposition \ref{keyaffine} we have to show that if $\gamma \in \posroots$
and $\al + \gamma \in \posroots$, then 
$k_{\al} + k_{\gamma} \leq k_{\al + \gamma} \leq k_{\al} + k_{\gamma} + 1$.
Since $k_{\al} = 1$, this simplifies to 
$k_{\gamma} + 1 \leq k_{\al + \gamma} \leq k_{\gamma} + 2$.
The second inequality is clearly a strict one.  For the first inequality,
write $\al + \gamma = 
\sum_{i=1}^{k+1} \gamma_i$ with $\gamma_i \in \posroots - \ideal$.
Then $\gamma = -\al + \sum_{i=1}^{k+1} \gamma_i$. 
Invoking Lemma \ref{reorder} and re-ordering the roots accordingly, 
we have $-\al + \gamma_1 \in \ro \cup \{ 0 \}$.  But $\al$ is maximal in 
$\posroots - \ideal$,
so $-\al + \gamma_1 \in \ro^- \cup \{ 0 \}$.  Continuing to use Lemma \ref{reorder},
we can write $\gamma = (-\al + \sum_{i=1}^{j} \gamma_i) + \sum_{i=j+1}^{k+1} \gamma_i$
where $j \geq 2$ and 
the expression in parentheses belongs to $\posroots - \ideal$.
It follows that $k_{\gamma} \leq k-1$ and so $k_{\gamma} \leq k_{\al + \gamma}-1$ 
as desired.

It follows that there exists $y \in W_a$ such that $y = w s_i$ 
where $s_i \in S$.  Moreover, $w(\aff{\al}_i) = \al + \delta$.
Hence $x(\al_i) = \al$ and $\al_i \in  B_{h-1, \lambda}$.

\end{proof}

The proof yields the following corollary.

\begin{cor} \label{goes_to_simple}
Let $\ideal$ be an ideal.   
\begin{enumerate}
\item 
Let $w=w_{min}$ be the minimal element of $\ideal$. 
If $\al \in \ideal$ is minimal, then $w^{-1}(\al + \delta)$
is a negative affine simple root.
\item 
If $\ideal$ is strictly positive and $w=w_{max}$ is the maximal element of $\ideal$, 
then $\al \in \posroots - \ideal$ is maximal 
implies that
$w^{-1}(\al + \delta)$ is an affine simple root.
\end{enumerate}
\end{cor}

We can now prove 

\begin{thm} \label{minimals}
Let $\ideal$ be an ideal.
\begin{enumerate}
\item Let $\ideal_{min}$ be the minimal elements of $\ideal$.
Then there exists $J \subset \Pi$ and $y \in W$ such that
$y(\ideal_{min}) = J$.  
\item Let $\ideal$ be strictly positive and let
$\ideal_{max}$ be the maximal elements of $\posroots - \ideal$.
Then there exists $J \subset \Pi$ and $y \in W$ such that
$y(\ideal_{max}) = J$.  
\end{enumerate}
\end{thm}

\begin{proof}
As in the proof of Theorem \ref{counting1}, let 
$\tilde{w} \in W_a$ be such that  
$\tilde{w}(\Delta^t) = t \bar{A}$ and 
hence so that $\tilde{w}(D^t) =  \colat \cap t \bar{A}$.
We note that if we write $\tilde{w} = \tilde{x} \tau_{\tilde{\lambda}}$
then $\tilde{x}(\ro_b \cup \ro_{b-h}) = \ph$ where
$\ph = \Pi \cup \{ -\theta \}$ \cite{fan}.

For the first part of the theorem: 
let $w$ be the minimal element of $ST(\ideal)$.  Write $w = x \tau_{\lambda}$ where $x \in W$
and $\lambda \in \colat$.
Then $x^{-1}(\ideal_{min}) = B_{h+1, \lambda}$ by Proposition \ref{simples}.  
Then $\tilde{x}(x^{-1}(\ideal_{min}))$ 
are exactly the elements of $\ph$ such that
$\tilde{w}(\lambda)$ (an element of $\colat$) lies on the corresponding wall of $(h+1) \bar{A}$.

In \cite{sommers:affineweylgp} it was proved that 
every such subset is conjugate by an element of $W$ 
to some $J \subset \Pi$.  We give a simpler proof here.

Write $\theta = \sum_{\al \in \Pi} c_{\al} \al$ and set $c_{-\theta}=1$. 
Given $J' \subset \ph$, 
let 
$$d = \text{ gcd }(c_{\al} \ | \ \al \in  \ph - J').$$
Recall (for example, from \cite{sommers:b-c}) that 
$J'$ is conjugate by an element of $W$ to a subset of $\Pi$ if and 
only if $d=1$.
Let $J' = \tilde{x}(x^{-1}(\ideal_{min}))$.  If $J' \subset \Pi$, there is nothing to prove.
Otherwise let $\mu = \tilde{w}(\lambda)$.
Since $-\theta \in J'$, we have $\langle \theta, \mu \rangle = h+1$.
Hence 
$$h+1 = \langle \sum_{\al \in \Pi} c_{\al} \al, \mu \rangle
=  \sum_{\al \in \Pi} c_{\al} \langle \al, \mu \rangle
=  \sum_{\al \in \Pi - J'} c_{\al} \langle \al, \mu \rangle$$
where the last equality holds since $\langle \al, \mu \rangle = 0$ if $\al \in \Pi \cap J'$.
Hence we see that $d$ divides $h+1$.  But 
it is known that each $c_{\al}$ divides $h$, hence $d$ divides $h$.  Thus
$d$ divides $h+1$ and $h$ and so $d=1$.
Note that the proof remains valid for any $t$ which is good, i.e. prime to 
all of the $c_{\al}$'s, by replacing $h+1$ above with such a $t$.

This completes the proof for the first part.  The second part follows 
in an analogous fashion.
\end{proof}

Given a $\al \in \ro$, let $e_{\al}$ be a non-zero element in the corresponding root space
in $\g$.  The following corollary is immediate from the theorem.

\begin{cor}
The nilpotent element $\sum_{\alpha \in \ideal_{min}} e_{\alpha}$ is regular
in a Levi subalgebra of $\g$.
\end{cor}

Let $J \subset \Pi$. 
We can use the results of \cite{sommers:affineweylgp} to enumerate the number of 
ideals $\ideal$ such that $\ideal_{min}$ is conjugate by some element of $W$ to $J$.

Let $W_J$ be the Weyl group generated by the simple reflections corresponding to 
the elements of $J$.
Let $V^{J}$ denote the subspace of $V$ fixed point-wise by all the elements of $W_J$.
Consider the set of hyperplanes 
$$\{ V^{J} \cap H_{\al, 0} \ | \ \al \in \posroots \text { and }
V^{J} \not\subset H_{\al,0} \}.$$
This defines a hyperplane arrangement in $V^{J}$. 
It is known that this arrangement is free \cite{douglass}, \cite{broer:free}
and hence that its characteristic polynomial
$p^J(t)$ factors as $\prod_{i=1}^{n-j} (t - m^J_i)$, where $n$ is the rank of $G$, $j$
is the cardinality of $J$, and $m^J_i$ are positive integers \cite{terao}.  
These integers, called the Orlik-Solomon exponents, were computed in
\cite{orlik-solomon}.
An alternative way to compute the Orlik-Solomon exponents
was given in \cite{sommers:affineweylgp}, section 5.
In fact those results and the proof of the above theorem yield

\begin{prop} \label{numbers}
Let $J \subset \Pi$.  Let $N(W_J)$ denote the normalizer of $W_J$ in $W$.
Let $$\chi_J(t) = \frac{1}{[N(W_J): W_J]} p^J(t).$$
\begin{enumerate}
\item The number of ideals $\ideal$ such that $\ideal_{min}$ is conjugate
under $W$ to $J$ is $\chi_J(h+1)$.
\item The number of strictly positive ideals $\ideal$ such that
$\ideal_{max}$ is conjugate
under $W$ to $J$ is $\chi_J(h-1)$.
\end{enumerate}
\end{prop}

\bibliography{sommers1}
\bibliographystyle{pnaplain}
\end{document}